\documentclass[12pt,reqno]{amsart}
\usepackage{amsmath}
\usepackage{graphicx}
\usepackage{hyperref}

\usepackage{mathtools}
\usepackage{amsfonts}
\usepackage{amssymb}
\usepackage[T1]{fontenc}
\usepackage{amsthm}
\usepackage{fullpage}
\usepackage{enumitem}
\usepackage{bbm}
\usepackage[usestackEOL]{stackengine}
\usepackage{multirow}
\newtheorem{theorem}{Theorem}[section]
\newtheorem{corollary}{Corollary}[section]
\newtheorem{lemma}[theorem]{Lemma}

\newtheorem{example}{Example}[section]
\newtheorem{proposition}{Proposition}[section]
\theoremstyle{definition}

\theoremstyle{remark}

\numberwithin{equation}{section}

\title{SCATTERED BEHAVIOR USING MODIFIED CYCLOTOMIC MAPPING OVER FINITE FIELDS OF ODD CHARACTERISTIC}


\author{Suman Mondal}
\address{Department of Mathematical Sciences, Tezpur University, Tezpur, Assam, 784028, India}
\email{mondalmondalsuman@gmail.com}

\begin{document}
	\begin{abstract}
		Introduced by Sheekey in 2016, the study of scattered polynomials over a finite field $\mathbb{F}_{q^n}$ has been increasing regarding the classification of those that are \textit{exceptional}, i.e., polynomials which are scattered over infinite field extensions, are limited to the cases where their index $t$ is small, or a prime number larger than the q-degree k of the polynomial, or an integer smaller than k in the case where k is a prime. In this paper, we focus on the scattered behavior of $S(x)=\sum_{i=1}^{{k}} a_ix^{q^{r_i}} \in \mathbb{F}_{q^n}[x]$, where $q$ is a power of an odd prime, $0<r_1<r_2< \cdots<r_k<n$ and $a_1, \cdots,a_k \in \mathbb{F}_{q^n}^*$ such that the order of $a_i$'s divide $(q^{r_1}-1)$, $\forall i=2,3,\cdots,k $. We explore a connection between $S(x)$ and 
        the cyclotomic mapping polynomial. As an application, in three parts, we discuss the scattered behavior of $S(x)$ of index $t$ where $t=r_1$, or $0<t<r_1$, or $r_1<t<n$. Starting with the pseudoregulus
        type of index $t \geq 0$, we present conditions to verify scattered behavior of $S(x)$ of index $r_1$. With some additional conditions, we do the same in case $0<t<r_1$ or $r_1<t<n$. In particular, for $S(x)=a_1x^{q^{r_1}}+a_2x^{q^{r_2}} \in \mathbb{F}_{q^n}[x]$ with $a_1,a_2 \in \mathbb{F}_{q^n}^*$ such that $|a_2| \mid q^{r_1}-1$, we present a necessary and sufficient condition to verify its scattered behavior of index $t \in \{r_1,r_2\}$. We also connect such scattered binomials with the well known \textit{Lunardon-Polverino} polynomial. With conditions on $\delta, q,n$, and $r$; we present a new family of \textit{exceptional scattered polynomial} $S(x)=x^q+\delta x^{q^{(2r+1)}} \in \mathbb{F}_{q^n}[x]$ of index $\{r+1\}$. 
	\end{abstract}
	
	\maketitle

\section{Introduction}
Let $p$ be a prime, $m\in \mathbb{Z}^+$, $q=p^m$, $n\in \mathbb{Z}^+$, and $S(x)=\sum_{i=0}^{{k}} a_ix^{q^i} \in \mathbb{F}_{q^n}[x]$ be an $\mathbb{F}_q$-linearized polynomial over $\mathbb{F}_{q^n}$. We also assume that the $q$-degree $k$ of $S(x)$ is smaller than n, so that the identification with the map $x \mapsto S(x)$ defines a one-to-one correspondence between such polynomials and $\mathbb{F}_{q}$-linear maps over $\mathbb{F}_{q^n}$. An $\mathbb{F}_{q}$-linearized polynomial $S(x) \in \mathbb{F}_{q^n}[x]$ is said to be a \textit{scattered polynomial} (SP)\cite{BARTOLI2018507} of index $t \in \{0,\cdots, n-1\} $, if for any distinct $y,z \in \mathbb{F}_{q^n}^*$,
    \begin{equation}\tag{1}\label{1}
        \frac{S(y)}{y^{q^{t}}} = \frac{S(z)}{z^{q^{t}}} \Rightarrow \frac{y}{z} \in \mathbb{F} _q.  
    \end{equation}
    Starting from \cite{BARTOLI2018507}, a much stronger property regarding the scattered polynomials, namely their \textit{exceptionality}, has been defined and deeply investigated. An $\mathbb{F}_q$-linearized polynomial
$S(x) \in  \mathbb{F}_{q^n}[x]$ is said to be \textit{exceptional scattered} of index $t \in \{0, \cdots, n -1\}$ if there exist
infinitely many $m \in \mathbb{N}$ such that, for any distinct $y,z \in  \mathbb{F}_{q^n}^*$, Condition (\ref{1}) holds. \\
 Scattered polynomials $S(x) \in \mathbb{F}_{q^n}[x]$ yield \textit{scattered subspaces} $U_S$ (w.r.t. a Desarguesian spread) in $\mathbb{F}_{q^n} \times \mathbb{F}_{q^n}$ by defining $$U_S= \{(x^{q^t}, S(x))| x \in \mathbb{F}_{q^n} \}. $$ 
 Scattered behavior of such polynomials have many applications, such as translation caps in affine spaces \cite{10.1007/s00493-016-3531-6}, blocking sets \cite{10.1006/ffta.2000.0280}, translation spreads of the Cayley generalized hexagon \cite{marino2015translation},
coding theory (\cite{polverino2020connections}, \cite{sheekey2015new}), and graph theory \cite{https://doi.org/10.1112/blms/18.2.97} etc.\\
Suppose $q$ is a power of an odd prime, $0<r_1<r_2< \cdots<r_k<n$ and $a_1, \cdots,a_k \in \mathbb{F}_{q^n}^*$ such that order of $a_i$'s divide $(q^{r_1}-1)$, $\forall i=2,3,\cdots,k $. We consider the  $\mathbb{F}_q$-linearized polynomial $S(x)=\sum_{i=1}^{{k}} a_ix^{q^{r_i}} \in \mathbb{F}_{q^n}[x]$. Then
\begin{align*}
         S(x)&=x^{q^{r_1}}[a_1+a_2x^{(q^{r_2}-q^{r_1})}+ \cdots + a_kx^{(q^{r_k}-q^{r_1})}] \\
         & = x^{q^{r_1}}[a_1+a_2x^ {\{q^{r_1}{{(q^{{r_2}-{r_1}}-1)}\}}}+ \cdots + a_kx^ {\{q^{r_1}{{(q^{{r_k}-{r_1}}-1)}\}}}] \\
         & = x^{q^{r_1}} f(x^{s \cdot q^{r_1}}), \text{where}\; s=q^i-1 \mid q^n-1, \text{for some}\; i \in \mathbb{N}.       
    \end{align*}
     With the conditions above, the $\mathbb{F}_q$-linearized polynomial $S(x)=\sum_{i=1}^{{k}} a_ix^{q^{r_i}} \in \mathbb{F}_{q^n}[x]$ is of the form $S(x)=  x^{q^{r_1}} f(x^{s \cdot q^{r_1}})$ where $s$ is a positive integer such that $s \mid q^n-1$. In this paper, we investigate the scattered behavior of such $S(x)\in \mathbb{F}_{q^n}[x]$, where $f(x)$ is an arbitrary polynomial of positive degree over $\mathbb{F}_{q^n}$ and $q^n-1=sl$, for some positive integer $l$ and $s$.
     
     We introduce the \textit{modified $r$-th order cyclotomic mapping} $f^{r}_{A_{0},\;A_{1},\;A_{2},\cdots, \;A_{l-1}}$ (derived from \cite{wang2007cyclotomic}), and establish an useful connection between polynomials of the form $x^{q^{r_1}} f(x^{s \cdot q^{r_1}})$ and the \textit{modified $r$-th order cyclotomic mapping} $f^{r}_{A_{0},\;A_{1},\;A_{2},\cdots, \;A_{l-1}}$ in Lemma(\ref{relation}). Several tools have already been proposed in the study of scattered polynomials and exceptional scattered polynomials, that are  
    related to certain algebraic curves or Galois extensions of function fields. Also, we have been investigating the exceptional scatteredness of a polynomial by considering separately the
    exceptionality of two weaker properties defined in \cite{longobardi2021partially}, namely the \textit{L-$q^t$-partial scatteredness} and the \textit{ R-$q^t$-partial scatteredness}.

    In this paper, in Section(\ref{sec2}), we start with the relation of polynomials of the form $S(x)=\sum_{i=1}^{{k}} a_ix^{q^{r_i}}=x^{q^{r_1}} f(x^{s \cdot q^{r_1}})$ and \textit{modified $r$-th order cyclotomic mapping} $f^{r}_{A_{0},\;A_{1},\;A_{2},\cdots, \;A_{l-1}}$. For $S(x)=\sum_{i=1}^{{k}} a_ix^{q^{r_i}} \in \mathbb{F}_{q^n}[x]$, we divide the index $t \in \{0, \cdots, n-1\}$ this into three parts, which are $t=r_1$, $0<t<r_1$, and $r_1<t<n$. 
    
    In Section(\ref{sec3}), we discuss a necessary and sufficient condition for $S(x)=\sum_{i=1}^{{k}} a_ix^{q^{r_i}} \in \mathbb{F}_{q^n}[x]$ to be a SP of index $r_1$. With some additional conditions, we do the same in case $0<t<r_1$, and $r_1<t<n$. From \cite{lunardon2014maximum}, we know that for $(r,n)=1$, $S(x)=x^{q^r}$(pseudoregulus type) is a SP of index $t=0$. In  Section(\ref{sec3}), we discuss the scattered behavior of $S(x)=x^{q^r}$, in case the index $t \geq 0$. From  Proposition (2.6) in \cite{longobardi2021partially}, we connect the scattered behavior of  $S(x)$ of index $t=0$ with  permutation behavior of a specific polynomial over $\mathbb{F}_{q^n}$. Here using the same proposition with the help of \textit{modified $r$-th order cyclotomic mapping}, we discuss the scattered behavior of $S(x)$ of index $t$, in terms of permutation behavior of a specific polynomial over $\mathbb{F}_{q^n}$, where $0<t<r_1$.

    In Section(\ref{sec4}), we consider the binomials of the form $S(x)=a_1x^{q^{r_1}}+a_2x^{q^{r_2}} \in \mathbb{F}_{q^n}[x]$ such that $|a_2| \mid q^{r_1}-1$. We discuss a necessary and sufficient condition for $S(x)$ to be a SP of index $t \in \{r_1,r_2\}$. Observe that $S(x)=a_1x^{q^{r_1}}+a_2x^{q^{r_2}} \in \mathbb{F}_{q^n}[x]$ and the Lunerdon-Polverino polynomials are not related in general, however, with some additional condition, in Lemma(\ref{transf}) we show that $S(x)$ can be transformed into an LP type polynomial. Also, from \cite{BARTOLI2022101956}, for $N_{q^n/q}(\delta) \neq 1$ and $(r,n)=1$,  $T(x)=x+\delta x^{q^{2r}} \in \mathbb{F}_{q^n}[x]$ is a exceptional scattered polynomial of index $\{r\}$. Using this result and the scattered properties that we obtain from \textit{modified $r$-th order cyclotomic mapping}, we work on a new family of \textit{exceptional scattered polynomial} in Theorem(\ref{exp}).
      \section{Modified Cyclotomic mapping}\label{sec2}
    Let $\gamma$ be a primitive element of  $(\mathbb{F}_{q^n}^\ast, \cdot)$ and $C_0$ be the collection of all the $l$-th power of $\gamma$, i.e., $C_0=\{ \gamma^0, \gamma^{1l}, \gamma^{2l}, \cdots ,\gamma^{(s-1)l}, \gamma^{sl}, \gamma^{(s+1)l}, \cdots  $ \} where $q^n-1=ls$ for some $l$, $s$ $\in$ $\mathbb{Z^+}$.
     Then,
    \begin{center}
    $C_0=\{ \gamma^{lj} : j=0, 1, 2, \cdots, s-1$\}.   
    \end{center}
    
    We observe that $ C_0 $ is a subgroup of the cyclic group $ ( \mathbb{F}_{q^n}^\ast , \cdot )$, so the quotient group $\mathbb{F}_{q^n}^\ast / C_0$ exists with respect to multiplication, with index $l$. The elements of $ \mathbb{F}_{q^n}^\ast / C_0$ are called the cyclotomic cosets $C_i$ and are defined by
    \begin{center}
    {$C_i=\gamma^i C_0,\;\; i=0, 1, 2, \cdots, l-1.$}
    \end{center}
    
    Let $x \in C_i$, then $x$ is of the form $\gamma^{i+lj}$ where $0 \leq i \leq l-1$ and $0 \leq j \leq s-1$. Suppose $r(=r_1)\in \mathbb{Z^+}$ and $ {A_0}, A_1, A_2, \cdots, A_{l-1} \in \mathbb{F}_{q^n}$, then we define \textit{modified $r$-th order cyclotomic mapping} $f^{r}_{A_{0},\;A_{1},\;A_{2},\cdots, \;A_{l-1}}$ from $\mathbb{F}_{q^n}$ to itself, as \\
   $$ f^{r}_{{A_0}, A_1, A_2, \cdots , A_{l-1}}(x)
   =\begin{cases}
       0 &  \text{for}  \ x=0\\
       A_i x^{q^{r_1}} & \text{for} \ x \in C_i.
   \end{cases}
   $$
    
    The polynomial $f^{r}_{A_{0},\;A_{1},\;A_{2},\cdots, \;A_{l-1}}(x)$ over $\mathbb{F}_{q^n}$, representing the \textit{modified $r$-th order cyclotomic mapping} $f^{r}_{A_{0},\;A_{1},\;A_{2},\cdots, \;A_{l-1}}$, is called a\textit{ modified $r$-th order cyclotomic mapping polynomial}. In particular, if $r=1$,  the obtained polynomial is known as \textit{modified cyclotomic mapping polynomial}.\\
    Let $\xi=\gamma^s$, then $\xi$ is a primitive $l$-th roots of unity. Now for $A_i \in \mathbb{F}_{q^n}$, we define $A_i=f (  \xi^{i \cdot q^{r_1}} )$ for $0\leq i\leq l-1$, where $\xi$ is a primitive $l$-th roots of unity.
    \begin{lemma}\label{relation}
        For any $r(=r_1)\in \mathbb{N}$, we have $S(x)=$ $x^{q^{r_1}}f(x^{s \cdot q^{r_1}})=f^{r}_{A_{0},\;A_{1},\;A_{2},\cdots, \;A_{l-1}}(x)$ where $A_i=f (  \xi^{i \cdot {q^{r_1}}} )$ for $0\leq i\leq l-1$,  $\xi$ is a primitive $l$-th roots of unity.
    \end{lemma}    
    \begin{proof}
     For $x=0$, we have $x^{q^{r_1}}f(x^{s \cdot q^{r_1}})=0=f^{r}_{A_{0},\;A_{1},\;A_{2},\;\cdots\;, \;A_{l-1}}(x)$, so the equality holds trivially.\\
     For $x  \in \mathbb{F}_{q^n}^\ast$, we have $x\in C_i$,   that is, $x$ is of the form $\gamma^{i+lj}$ for $0\leq i\leq l-1$ and $0\leq j\leq s-1$.
    Now, $f(x^{s \cdot q^{r_1}})=f(\gamma^{s \cdot q^{r_1} \cdot (i+lj)})=f(\gamma^{isq^{r_1}})=f(\xi^{i \cdot q^{r_1}})=A_i$, for $0\leq i\leq l-1$.\\
    Hence, $x^{q^{r_1}}f(x^{s \cdot q^{r_1}})=f^{r}_{A_{0},\;A_{1},\;A_{2},\cdots, \;A_{l-1}}(x)$ where $A_i=f (  \xi^{i \cdot {q^{r_1}}} )$ for $0\leq i\leq l-1$,  $\xi$ is a primitive $l$-th roots of unity.
    \end{proof}
    \begin{example}
        Let $S(x)=x^{9^2}+x^{9^4}+x^{9^6} \in \mathbb{F}_{9^8}[x]$. Then\\
        $S(x)=x^{9^2}+x^{9^4}+x^{9^6}=x^{9^2}[1+x^{9^2 \cdot (9^2-1)}+x^{9^2 \cdot (9^4-1)}]=x^{9^2} f(x^{9^2 \cdot (9^2-1)})$. For $f(x)=1+x+x^{(9^2+1)}$, $s=9^2-1$, and $l=\frac{9^8-1}{9^2-1}$, we have $S(x)=x^{9^{2}}f(x^{(9^2-1) \cdot 9^{2}})=f^{2}_{A_{0},\;A_{1},\;A_{2},\cdots, \;A_{l-1}}(x)$.
    \end{example}
    \begin{example}
        Let $S(x)=x^{5^3}+x^{5^4} \in \mathbb{F}_{5^5}[x]$. Then
        $S(x)=x^{5^3}+x^{5^4}=x^{5^3}[1+x^{5^3 \cdot (5-1)}]=x^{5^3} f(x^{5^3 \cdot (5-1)})$. For $f(x)=1+x$, $s=5-1$, and $l=\frac{5^5-1}{5-1}$, we have $S(x)=x^{5^{3}}f(x^{(5-1) \cdot 5^{3}})=f^{3}_{A_{0},\;A_{1},\;A_{2},\cdots, \;A_{l-1}}(x)$.
    \end{example}
    Suppose $\gamma$ is a primitive element of  $(\mathbb{F}_{q^n}^\ast\cdot)$, then for $x (\neq 0) \in \mathbb{F}_{q^n}$, $x$ is of the form $x=\gamma^a$, for some $a \in \mathbb{N}$. If $y \in \mathbb{F}_q$, then trivially $y \in \mathbb{F}_{q^n}$. However, for $\gamma^a \in \mathbb{F}_{q^n}$, we do not know if $\gamma^a \in \mathbb{F}_q$ or $\gamma^a \in \mathbb{F}_{q^n} - \mathbb{F}_q$. Next we present a result where for $\gamma^a \in \mathbb{F}_{q^n}$ we have a condition involving $a,q$ and $n$ such that $\gamma^a \in \mathbb{F}_q$.
    \begin{lemma}\label{power}
        Consider the extension $\mathbb{F}_{q^n} / \mathbb{F}_{q}$ and let $\gamma$ be a primitive element of  $(\mathbb{F}_{q^n}^\ast, \cdot)$ with $a\in \mathbb{N}$. Then $\gamma^a \in \mathbb{F}_{q}$ iff $\displaystyle \sum_{i=0}^{n-1} q^i \mid a$.
    \end{lemma}
    \begin{proof}
        Let $\gamma^a \in \mathbb{F}_{q}$, then $\gamma^{a(q-1)}=1$, i.e, $\displaystyle \sum_{i=0}^{n-1} q^i \mid a$.

        Conversely, let $\displaystyle \sum_{i=0}^{n-1} q^i \mid a$. Then $q^n-1 \mid a(q-1)$.\\
        Now $\mathbb{F}_{q^n} / \mathbb{F}_{q}$ has a normal basis of the form $\{ \gamma, \gamma^q, \cdots, \gamma^{q^{n-1}}\}$. So, for $\gamma^a \in \mathbb{F}_{q^n}$, we have $\gamma^a= \displaystyle \sum_{i=0}^{n-1} C_i\gamma^{q^i}$ where $C_i \in \mathbb{F}_{q}$. Also, $C_i^q=C_i$.\\
        Now
            $\gamma^a= \displaystyle \sum_{i=0}^{n-1} C_i\gamma^{q^i}$
             $\Rightarrow  \gamma^{a(q-1)}=1= (\displaystyle \sum_{i=0}^{n-1} C_i\gamma^{q^i})^{q-1}$
             $\Rightarrow  (\displaystyle \sum_{i=0}^{n-1} C_i\gamma^{q^i})^{q}=\displaystyle \sum_{i=0}^{n-1} C_i\gamma^{q^i}$
              $\Rightarrow  $ $\displaystyle \sum_{i=0}^{n-1} C_i^q\gamma^{q^{i+1}}$ $=\displaystyle \sum_{i=0}^{n-1} C_i\gamma^{q^i}$
             $\Rightarrow  $ $(\displaystyle \sum_{i=0}^{n-1} C_i\gamma^{q^{i+1}})=\displaystyle \sum_{i=0}^{n-1} C_i\gamma^{q^i}$
             $\Rightarrow C_0 \cdot\gamma+C_1\cdot\gamma^q+C_2\cdot\gamma^{q^2}+ \cdots+C_{n-1}\cdot\gamma^{q^{n-1}}=C_0 \cdot\gamma^q+C_1\cdot\gamma^{q^2}+C_2\cdot\gamma^{q^3}+ \cdots+C_{n-2}\cdot\gamma^{q^{n-1}}+C_{n-1}\cdot\gamma$
             $\Rightarrow (C_0-C_{n-1}) \cdot\gamma+(C_1-C_0)\cdot\gamma^q+(C_2-C_1)\cdot\gamma^{q^2}+ \cdots+(C_{n-1}-C_{n-2})\cdot\gamma^{q^{n-1}}=0$
             $\Rightarrow C_0-C_{n-1}= 0 = C_i-C_{i-1}, \forall i=1,2, \cdots, n-1.$\\
             Let $C_0=C_1= \cdots= C_{n-1}= C \in \mathbb{F}_q$.
             Then $\gamma^a= C\displaystyle \sum_{i=0}^{n-1} \gamma^{q^i}=C \cdot Tr_{(q^n/q)}(\gamma) \in \mathbb{F}_q$.\\
        So, $\gamma^a \in \mathbb{F}_q$.
    \end{proof}
    \section{Scattered behavior using cyclotomic mapping}\label{sec3}
    Suppose $S(x)=\sum_{i=1}^{{k}} a_ix^{q^{r_i}} \in \mathbb{F}_{q^n}[x]$ is a SP of index $t \in \{0, \cdots, n-1\}$. Then for some distinct $y,z \in \mathbb{F}_{q^n}^*$, $(y,z)$ satisfies (\ref{1}). From now on, by deciding pair (DP), we mean there exists $(y,z) \in C_i \times C_j \subset \mathbb{F}_{q^n}^* \times \mathbb{F}_{q^n}^*$ such that $(y,z)$ satisfies (\ref{1}). As every linearized polynomial of the form $S(x)=\sum_{i=1}^{{k}} a_ix^{q^{r_i}} \in \mathbb{F}_{q^n}[x]$ can be written as $S(x)=\sum_{i=1}^{{k}} a_ix^{q^{r_i}}=x^{q^r}A_i=f^{r}_{A_{0},\;A_{1},\;A_{2},\cdots, \;A_{l-1}}(x) \in \mathbb{F}_{q^n}[x]$, then for some $y(\in C_i),z(\in C_j) \in \mathbb{F}_{q^n}^*$, we have
    \begin{equation}\tag{2}\label{2}
        \frac{y^{q^r}A_i}{y^{q^{t}}} = \frac{z^{q^r}A_j}{z^{q^{t}}} \Rightarrow \frac{y}{z} \in \mathbb{F} _q.  
    \end{equation}
    For a given SP of any index over $\mathbb{F}_{q^n}$, and for a given DP $(y,z) \in C_i \times C_j \subset \mathbb{F}_{q^n}^* \times \mathbb{F}_{q^n}^*$, From (\ref{2}) we try to obtain a relation between $A_i$ and $A_j$.
     \begin{theorem}\label{ai}
      Let $S(x)$ be a SP of index $t \in \{0, \cdots, n-1\}$ over $\mathbb{F}_{q^n}$ and $(y,z)$ be a DP for some $i,j$ with $0 \leq i,j \leq l-1$ such that $y \in C_i$ and $ z \in C_j$. Then, $A_i=A_j$.
  \end{theorem}
  \begin{proof}
      If $i=j$, then trivially $A_i=A_j$. Also, for $A_i=A_j=0$, we have $A_i=A_j$.\\
      Suppose $i \neq j$ and $A_i,A_j \in \mathbb{F}_{q^n}^*$. Then we prove $A_i=A_j$.\\
      If possible let $A_i \neq A_j$, i.e., $A_i/A_j \neq 1$.\\
      As $S(x)$ is a SP of index $t \in \{0, \cdots, n-1\}$ over $\mathbb{F}_{q^n}$, let $(y,z)$ be a DP (for some $i,j$ with $0 \leq i,j \leq l-1$ such that $y \in C_i$, $z \in C_j$, and $y= \gamma^{lj_1+i}, z=\gamma^{lj_2+j}$, for some $j_1,j_2$ with $0 \leq j_1,j_2 \leq s-1$) such that $\frac{S(y)}{y^{q^{t}}} = \frac{S(z)}{z^{q^{t}}}$. Now
      \begin{align*}
          & \frac{S(y)}{y^{q^{t}}} = \frac{S(z)}{z^{q^{t}}}\\
           \Rightarrow \; & \frac{y^{q^{r_1}}A_i}{y^{q^{t}}} = \frac{z^{q^{r_1}}A_j}{z^{q^{t}}}\\
           \Rightarrow \; & (\frac{y}{z})^{(q^{r_1}-q^t)}=A_i/A_j \neq 1\\
           \Rightarrow \; & \gamma^{[l(j_1-j_2)+(i-j)](q^{r_1}-q^t)} \neq 1\\
           \Rightarrow \; & q^n-1 \nmid [l(j_1-j_2)+(i-j)](q^{r_1}-q^t)\\
           \Rightarrow \; & \displaystyle \sum_{i=0}^{n-1} q^i \nmid l(j_1-j_2)+(i-j)\\
           \Rightarrow \; & \gamma^{l(j_1-j_2)+(i-j)}=\frac{y}{z} \not\in \mathbb{F}_q, \; \text{[using Lemma (\ref{power})]}
      \end{align*}
      which is not possible, as $S(x)$ is a SP of index $t \in \{0, \cdots, n-1\}$ over $\mathbb{F}_{q^n}$.\\
      Hence $A_i=A_j$.
  \end{proof}

   We know that the Psudoregulus type $S(x)=x^{q^r}$ is of index $0$ over $\mathbb{F}_{q^n}$ \textit{iff} $(r,n)=1$\cite{lunardon2014maximum}. Next, we discuss the scattered behavior that involves positive index (including zero) of $S(x)=x^{q^r}$ over $\mathbb{F}_{q^n}$.
    \begin{proposition}\label{pseudo}
    Let $r \in \mathbb{N}$  and $t( \neq r)$ be a non negative integer, then $S(x)=x^{q^r}$ is a SP of index $t$ over $\mathbb{F}_{q^n}$ iff $(|t-r|, n)=1$. 
    \end{proposition}
    \begin{proof}
       We know that for $r \in \mathbb{N}$, $x^{q^r}$ is a SP of index $0$ over $\mathbb{F}_{q^n}$ \textit{iff} $(r, n)=1$. So, the proposition is true for $t=0$. \\
       Suppose $t(\neq r) >0$, and let $(y,z) \in C_i \times C_j \subset \mathbb{F}_{q^n}^* \times \mathbb{F}_{q^n}^*$ be a DP. Then
        \begin{align*}
          \frac{S(y)}{y^{q^{t}}} = \frac{S(z)}{z^{q^{t}}}
           \Leftrightarrow \;  y^{(q^{r_1}-q^t)}=z^{(q^{r_1}-q^t)}
           \Leftrightarrow \;  (\frac{y}{z})^{(q^{|t-r|}-1)}=1
           \Leftrightarrow \; \frac{y^{q^{|t-r|}}}{y}=\frac{z^{q^{|t-r|}}}{z} \tag{3}\label{3}
      \end{align*}
      Let  $S(x)=x^{q^r}$ be a SP of index $t$ over $\mathbb{F}_{q^n}$, then for some DP $(y_1, z_1) \in \mathbb{F}_{q^n}^* \times \mathbb{F}_{q^n}^*$, from $(\ref{1})$ we have  
       \begin{align*}
        \frac{S(y_1)}{y_1^{q^{t}}} = \frac{S(z_1)}{z_1^{q^{t}}} \Rightarrow \frac{y_1}{z_1} \in \mathbb{F} _q.  
    \end{align*}
    So, from (\ref{3}) we have $$\frac{y_1^{q^{|t-r|}}}{y_1}=\frac{z_1^{q^{|t-r|}}}{z_1} \Rightarrow \frac{y_1}{z_1} \in \mathbb{F} _q.$$
Hence, $x^{q^{|t-r|}}$ is a SP of index $0$, i.e., $(|t-r|, n)=1$.\\
    Conversely, let $(|t-r|, n)=1$, i.e., $x^{q^{|t-r|}}$ is a SP of index $0$ over $\mathbb{F}_{q^n}$. Then for some DP $(y_2, z_2) \in \mathbb{F}_{q^n}^* \times \mathbb{F}_{q^n}^*$, from (\ref{1}) we have $$\frac{y_2^{q^{|t-r|}}}{y_2}=\frac{z_2^{q^{|t-r|}}}{z_2} \Rightarrow \frac{y_2}{z_2} \in \mathbb{F} _q$$.
    So, from (\ref{3}), $S(x)=x^{q^r}$ is a SP of index $t$ over $\mathbb{F}_{q^n}$.
    \end{proof}
    \begin{example}
    $S(x)=x^{25^8} \in \mathbb{F}_{25^{15}}[x]$ is a SP of index $t \in \{ 0,1,4,6,7,9,10,12 \}$.
    \end{example}
     Next we discuss the scattered behavior for index $t=r_1$ of $S(x) \in \mathbb{F}_{q^n}[x]$. Also, we know that if $S(x)$ be a SP of any index, then there exists a DP satisfying (\ref{1}). However, out of many such pairs, explicitly we do not know if that DP is unique or not. In the next result, we observe that, in some cases, the DP may not be unique.
    \begin{theorem}\label{1.5}
      Let $S(x)=\sum_{i=1}^{{k}} a_ix^{q^{r_i}} \in \mathbb{F}_{q^n}[x]$ such that $|a_i| \mid q^{r_1}-1$, $\forall i=2, \cdots, k$. Then the following are equivalent.
      \begin{enumerate}
          \item  $S(x)=\sum_{i=1}^{{k}} a_ix^{q^{r_i}} \in \mathbb{F}_{q^n}[x]$  is a SP of index $r_1$,
          \item $S_{r_1}(x)=S(x)-a_1x^{q^{r_1}}=\sum_{i=2}^{{k}} a_ix^{q^{r_i}} \in \mathbb{F}_{q^n}[x]$  is a SP of index $r_1$,
          \item $S_{r_1}^{r_1} (x)=\sum_{i=2}^{{k}} a_ix^{q^{(r_i-r_1)}} \in \mathbb{F}_{q^n}[x]$ is a SP of index $0$.
      \end{enumerate}
      Further, suppose $n$ is odd, $(r,n)=1$, and $S(x)=\sum_{i=1}^{{k}} a_ix^{q^{r_i}} \in \mathbb{F}_{q^n}[x]$  is a SP of index $r_1$, then DP may not be unique.
    \end{theorem}
    \begin{proof}
    Let $(y,z) \in C_i \times C_j \subset \mathbb{F}_{q^n}^* \times \mathbb{F}_{q^n}^*$. Then
        \begin{align*}
          & \frac{S(y)}{y^{q^{r_1}}} = \frac{S(z)}{z^{q^{r_1}}}
           \Leftrightarrow \; \frac{\sum_{i=2}^{{k}} a_iy^{q^{r_i}}}{y^{q^{r_1}}}=\frac{\sum_{i=2}^{{k}} a_iz^{q^{r_i}}}{z^{q^{r_1}}}
           \Leftrightarrow \; \frac{S_{r_1}(y)}{y^{q^{r_1}}} = \frac{S_{r_1}(z)}{z^{q^{r_1}}}. \tag{4} \label{4}
      \end{align*}
        Let $S(x)=\sum_{i=1}^{{k}} a_ix^{q^{r_i}} \in \mathbb{F}_{q^n}[x]$  be a SP of index $r_1$ with DP $(y_1,z_1)$. Then from (\ref{1}) we have 
        $$\frac{S(y_1)}{y_1^{q^{r_1}}} = \frac{S(z_1)}{z_1^{q^{r_1}}} \Rightarrow \frac{y_1}{z_1} \in \mathbb{F}_q. $$ 
        Now from (\ref{4}) we have $$\frac{S_{r_1}(y_1)}{y_1^{q^{r_1}}} = \frac{S_{r_1}(z_1)}{z_1^{q^{r_1}}} \Rightarrow \frac{y_1}{z_1} \in \mathbb{F} _q. $$
        So, $S_{r_1}(x)=\sum_{i=2}^{{k}} a_ix^{q^{r_i}} \in \mathbb{F}_{q^n}[x]$  is a SP of index $r_1$. Conversely, if $S_{r_1}(x)=\sum_{i=2}^{{k}} a_ix^{q^{r_i}} \in \mathbb{F}_{q^n}[x]$  is a SP of index $r_1$, from (\ref{4}), it implies that $S(x)=\sum_{i=1}^{{k}} a_ix^{q^{r_i}} \in \mathbb{F}_{q^n}[x]$  is a SP of index $r_1$.\\
        Hence, (1) and (2) are equivalent.\\
    Suppose $(Y,Z) \in  \mathbb{F}_{q^n}^* \times \mathbb{F}_{q^n}^*$. Then
    \begin{align*}
            & \frac{S_{r_1}(Y)}{Y^{q^{r_1}}} = \frac{S_{r_1}(Z)}{Z^{q^{r_1}}}\\
            \Leftrightarrow \; &  \frac{a_2Y^{q^{r_2}}+ \cdots +a_kY^{q^{r_k}}}{Y^{q^{r_1}}} = \frac{a_2Z^{q^{r_2}}+ \cdots +a_kZ^{q^{r_k}}}{Z^{q^{r_1}}}\\
            \Leftrightarrow \; & a_2Y^{(q^{r_2}-q^{r_1})}+ \cdots +a_kY^{(q^{r_k}-q^{r_1})}=a_2Z^{(q^{r_2}-q^{r_1})}+ \cdots +a_kZ^{(q^{r_k}-q^{r_1})}\\
           \Leftrightarrow \; & {a_2Y^ {\{q^{r_1}{{(q^{{r_2}-{r_1}}-1)}\}}}+ \cdots + a_kY^ {\{q^{r_1}{{(q^{{r_k}-{r_1}}-1)}\}}}}={a_2Z^ {\{q^{r_1}{{(q^{{r_2}-{r_1}}-1)}\}}}+ \cdots + a_kZ^ {\{q^{r_1}{{(q^{{r_k}-{r_1}}-1)}\}}}}\\
           \Leftrightarrow \; & {a_2Y^ {{{(q^{{r_2}-{r_1}}-1)}}}+ \cdots + a_kY^ {{{(q^{{r_k}-{r_1}}-1)}}}} = {a_2Z^ {{{(q^{{r_2}-{r_1}}-1)}}}+ \cdots + a_kZ^ {{{(q^{{r_k}-{r_1}}-1)}}}}\\
           \Leftrightarrow \; & \frac{S_{r_1}^{r_1}(Y)}{Y} = \frac{S_{r_1}^{r_1}(Z)}{Z}.\tag{5} \label{5}
            \end{align*}
            Let $S_{r_1}(x)=\sum_{i=2}^{{k}} a_ix^{q^{r_i}} \in \mathbb{F}_{q^n}[x]$ be a SP of index $r_1$  with DP $(y_2,z_2)$. Then from (\ref{1}) we have $$\frac{S_{r_1}(y_2)}{y_2^{q^{r_1}}} = \frac{S_{r_1}(z_2)}{z_2^{q^{r_1}}} \Rightarrow \frac{y_2}{z_2} \in \mathbb{F} _q. $$
            Now from (\ref{5}) we have $$\frac{S_{r_1}^{r_1}(y_2)}{y_2} = \frac{S_{r_1}^{r_1}(z_2)}{z_2} \Rightarrow \frac{y_2}{z_2} \in \mathbb{F} _q.$$
        So, $S_{r_1}^{r_1}(x)=\sum_{i=2}^{{k}} a_ix^{q^{(r_i-r_1)}}$ $ \in \mathbb{F}_{q^n}[x]$  is a SP of index $0$. Conversely, if $S_{r_1}^{r_1}(x)=\sum_{i=2}^{{k}} a_ix^{q^{(r_i-r_1)}}$ $ \in \mathbb{F}_{q^n}[x]$  is a SP of index $0$, from (\ref{5}), it implies that $S_{r_1}(x)=\sum_{i=2}^{{k}} a_ix^{q^{r_i}} \in \mathbb{F}_{q^n}[x]$  is a SP of index $r_1$.\\
        Hence, (2) and (3) are equivalent, implying  (1), (2), (3) are equivalent.

            For the second part, we have $n$ is odd, and $(r,n)=1$. Suppose $S(x)=\sum_{i=1}^{{k}} a_ix^{q^{r_i}} \in \mathbb{F}_{q^n}[x]$  is a SP of index $r_1$ and let $(y,z) \in C_i \times C_j \subset \mathbb{F}_{q^n}^* \times \mathbb{F}_{q^n}^*$ be a DP for $S(x)$. Then from (\ref{4}), we have $$\frac{S_{r_1}(y)}{y^{q^{r_1}}} = \frac{S_{r_1}(z)}{z^{q^{r_1}}} \Rightarrow \frac{y}{z} \in \mathbb{F} _q. $$. So, $(y,z)$ is also a DP for $S_{r_1}(x)$ of index $r_1$. Let $Y=y^{q^{r_1}}, Z=z^{q^{r_1}}$. Then
    \begin{align*}
            & \frac{S(y)}{y^{q^{r_1}}} = \frac{S(z)}{z^{q^{r_1}}}\\
           \Leftrightarrow \; & {a_2y^ {\{q^{r_1}{{(q^{{r_2}-{r_1}}-1)}\}}}+ \cdots + a_ky^ {\{q^{r_1}{{(q^{{r_k}-{r_1}}-1)}\}}}}={a_2z^ {\{q^{r_1}{{(q^{{r_2}-{r_1}}-1)}\}}}+ \cdots + a_kz^ {\{q^{r_1}{{(q^{{r_k}-{r_1}}-1)}\}}}}\\
           \Leftrightarrow \; & \frac{a_2Y^ {{{q^{{r_2}-{r_1}}}}}+ \cdots + a_kY^ {{{q^{{r_k}-{r_1}}}}}}{Y} = \frac{a_2Z^ {{{q^{{r_2}-{r_1}}}}}+ \cdots + a_kZ^ {{{q^{{r_k}-{r_1}}}}}}{Z}\\
           \Leftrightarrow \; & \frac{S_{r_1}^{r_1}(Y)}{Y} = \frac{S_{r_1}^{r_1}(Z)}{Z}\\
           \Leftrightarrow \; & \frac{S_{r_1}(Y)}{Y^{q^{r_1}}} = \frac{S_{r_1}(Z)}{Z^{q^{r_1}}}[\text{using } (\ref{5})].\tag{6} \label{6}
            \end{align*} 
           As $S(x)=\sum_{i=1}^{{k}} a_ix^{q^{r_i}} \in \mathbb{F}_{q^n}[x]$  is a SP of index $r_1$ and $(y,z) \in C_i \times C_j \subset \mathbb{F}_{q^n}^* \times \mathbb{F}_{q^n}^*$ is a DP, then from (\ref{6}) we have 
           \begin{align*}
               \frac{S_{r_1}(Y)}{Y^{q^{r_1}}} = \frac{S_{r_1}(Z)}{Z^{q^{r_1}}} \Leftrightarrow \; \frac{S(y)}{y^{q^{r_1}}} = \frac{S(z)}{z^{q^{r_1}}} \Rightarrow \frac{y}{z} \in \mathbb{F}_q \Rightarrow \frac{Y}{Z} \in \mathbb{F}_q. 
           \end{align*}
           So, $(Y,Z)$ is also a DP for $S_{r_1}(x)$ of index $r_1$ and trivially $Y\neq Z$. If possible let $(y,z)=(Y,Z)$.\\
           Case 1: Let $y=Y=y^{q^{r_1}}$ and  $z=Z=z^{q^{r_1}}$. Then $ |y|,|z| \mid (q^{r_1}-1)$.\\
           So, $ |y|,|z| \mid (q^{r_1}-1,q^n-1)$.
           As $(r_1,n)=1$, then $(q^{r_1}-1,q^n-1)=q-1$, i.e., $ |y|,|z| \mid q-1$.\\
          Now from (\ref{5}), as $\frac{S_{r_1}^{r_1}(Y)}{Y} = \frac{S_{r_1}^{r_1}(Z)}{Z}$, then $1+\frac{a_2Y^ {{{q^{{r_2}-{r_1}}}}}+ \cdots + a_kY^ {{{q^{{r_k}-{r_1}}}}}}{Y} = 1+\frac{a_2Z^ {{{q^{{r_2}-{r_1}}}}}+ \cdots + a_kZ^ {{{q^{{r_k}-{r_1}}}}}}{Z}$, i.e., $Y=Z$, which is not possible.\\
          Case 2: Let $y=Z=z^{q^{r_1}}$ and  $z=Y=y^{q^{r_1}}$. Then $ |y|,|z| \mid (q^{2r_1}-1)$.\\
           So, $ |y|,|z| \mid (q^{2r_1}-1,q^n-1)$.
           As $(r_1,n)=1$ and $n$ is odd, then $(q^{2r_1}-1,q^n-1)=q-1$, i.e., $ |y|,|z| \mid q-1$. Similar as case 1, we obtain that 
           $Y=Z$, which is not possible.\\
           From both cases, we have $(y,z) \neq (Y,Z)$, where both $(y,z)$ and $(Y,Z)$ are DP for $S_{r_1}(x)$  of index $r_1$. So, DP may not be unique.
    \end{proof}
    Throughout this paper, for $S(x)=x^{q^{r_1}}f(x^{s \cdot q^{r_1}})$, we divide the index $t$ into 3 parts. One part where we consider $r_1=t$ and others are $0 \leq t <r_1$ and $r_1 < t <n$. In Theorem (\ref{1.5}), we considered the case $t=r_1$. Next, including some additional conditions, we discuss the other two cases.
 \begin{theorem}\label{22}
      Let $S(x)=\sum_{i=1}^{{k}} a_ix^{q^{r_i}} \in \mathbb{F}_{q^n}[x]$ such that $0<t<n$ and $|a_i| \mid q^{t}-1$ 
      (in case $0<t<r_1$), $|a_i| \mid q^{r_1}-1$ 
      (in case $r_1<t<n$) $\forall i=1, \cdots, k$. Then we have the following.
      \begin{enumerate}
          \item  $S(x)=\sum_{i=1}^{{k}} a_ix^{q^{r_i}} \in \mathbb{F}_{q^n}[x]$  is a SP of index $t(<r_1)$ \textit{iff} $S^t(x)=\sum_{i=1}^{{k}} a_ix^{q^{(r_i-t)}}$ is a SP of index 0,
          \item $S(x)=\sum_{i=1}^{{k}} a_ix^{q^{r_i}} \in \mathbb{F}_{q^n}[x]$  is a SP of index $t(>r_1)$ \textit{iff} $T(x)=a_1x+S_{r_1}^{r_1}(x)=a_1x+\sum_{i=2}^{{k}} a_ix^{q^{(r_i-r_1)}}$ is a SP of index $\{t-r_1\}$ over $\mathbb{F}_{q^n}$.
      \end{enumerate}
    \end{theorem}
    \begin{proof}
        For $t=0$, (1) is true trivially.\\
        Let $0<t<r_1$ and $(y,z) \in C_i \times C_j \subset \mathbb{F}_{q^n}^* \times \mathbb{F}_{q^n}^*$. Then
        \begin{align*}
          & \frac{S(y)}{y^{q^{t}}} = \frac{S(z)}{z^{q^{t}}}\\
           \Leftrightarrow \; & \sum_{i=1}^{{k}} a_iy^{(q^{r_i}-q^t)}=\sum_{i=1}^{{k}} a_iz^{(q^{r_i}-q^t)}\\
           \Leftrightarrow \; & \sum_{i=1}^{{k}} a_iy^{q^t(q^{r_i-t}-1)}=\sum_{i=1}^{{k}} a_iz^{q^t(q^{r_i-t}-1)}\\
           \Leftrightarrow \; & \frac{\sum_{i=1}^{{k}} a_iy^{q^{r_i-t}}}{y} = \frac{\sum_{i=1}^{{k}} a_iz^{q^{r_i-t}}}{z} [\text{as order of } a_i's \text{ divide } (q^t-1)]\\
           \Leftrightarrow \;  & \frac{S^t(y)}{y} = \frac{S^t(z)}{z}.\tag{7} \label{7}
      \end{align*}
      Let $S(x)$ be a SP of index $t$ with DP $(y_1,z_1)$. Then from (\ref{1}) we have 
        
        $$\frac{S(y_1)}{y_1^{q^{t}}} = \frac{S(z_1)}{z_1^{q^{t}}} \Rightarrow \frac{y_1}{z_1} \in \mathbb{F}_q. $$ 
        Now from (\ref{7}), we have $$\frac{S^t(y_1)}{y_1} = \frac{S^t(z_1)}{z_1} \Rightarrow  \frac{y_1}{z_1} \in \mathbb{F} _q. $$
        So, $S^t(x)=\sum_{i=1}^{{k}} a_ix^{q^{(r_i-t)}} \in \mathbb{F}_{q^n}[x]$  is a SP of index $0$. Conversely, if $S^t(x)=\sum_{i=1}^{{k}} a_ix^{q^{(r_i-t)}} \in \mathbb{F}_{q^n}[x]$  is a SP of index $0$, from (\ref{7}), it implies that $S(x)=\sum_{i=1}^{{k}} a_ix^{q^{r_i}} \in \mathbb{F}_{q^n}[x]$  is a SP of index $t$.\\
        Hence, (1) is proved.\\
        For the second part, let $r_1<t<n$ such that $t=r_1+u$, and $(y,z) \in C_i \times C_j \subset \mathbb{F}_{q^n}^* \times \mathbb{F}_{q^n}^*$. Then
        \begin{align*}
          & \frac{S(y)}{y^{q^{t}}} = \frac{S(z)}{z^{q^{t}}}\\
           \Leftrightarrow \; & \frac{a_1+\sum_{i=2}^{{k}} a_iy^{(q^{r_i}-q^{r_1})}}{y^{(q^{r_1+u}-q^{r_1})}} = \frac{a_1+\sum_{i=2}^{{k}} a_iz^{(q^{r_i}-q^{r_1})}}{z^{(q^{r_1+u}-q^{r_1})}}\\
           \Leftrightarrow \; & \frac{a_1+\sum_{i=2}^{{k}} a_iy^{(q^{r_i-r_1}-1)}}{y^{(q^u-1)}} = \frac{a_1+\sum_{i=2}^{{k}} a_iz^{(q^{r_i-r_1}-1)}}{z^{(q^u-1)}} [\text{as order of } a_i's \text{ divide } (q^{r_1}-1) ]\\
           \Leftrightarrow \; & \frac{a_1y+\sum_{i=2}^{{k}} a_iy^{q^{r_i-r_1}}}{y^{q^{t-r_1}}}=\frac{a_1y+\sum_{i=2}^{{k}} a_iz^{q^{r_i-r_1}}}{z^{q^{t-r_1}}}\\
           \Leftrightarrow \; & \frac{T(y)}{y^{q^{t-r_1}}}=\frac{T(z)}{z^{q^{t-r_1}}}.\tag{8}\label{8}
        \end{align*} 
        Let $S(x)$ be a SP of index $t(>r_1)$ with DP $(y_2,z_2)$. Then from (\ref{1}) we have 
        
        $$\frac{S(y_2)}{y_2^{q^{t}}} = \frac{S(z_2)}{z_2^{q^{t}}} \Rightarrow \frac{y_2}{z_2} \in \mathbb{F}_q. $$ 
        Now from (\ref{8}) we have $$\frac{T(y_2)}{y_2^{q^{t-r_1}}} = \frac{T(z_2)}{z_2^{q^{t-r_1}}} \Rightarrow \frac{y_2}{z_2} \in \mathbb{F} _q. $$
        So, $T(x)=a_1x+\sum_{i=2}^{{k}} a_ix^{q^{(r_i-r_1)}} \in \mathbb{F}_{q^n}[x]$  is a SP of index $\{t-r_1\}$. Conversely, if $T(x)=a_1x+\sum_{i=2}^{{k}} a_ix^{q^{(r_i-r_1)}} \in \mathbb{F}_{q^n}[x]$  is a SP of index $\{t-r_1\}$, from (\ref{8}), it implies that $S(x)=\sum_{i=1}^{{k}} a_ix^{q^{r_i}} \in \mathbb{F}_{q^n}[x]$  is a SP of index $t(>r_1)$.\\
        Hence, (2) is proved.
    \end{proof}
    In Proposition 2.6 of  \cite{longobardi2021partially}, for $t=0$, we connected the scatter behavior of index $t$ with the permutation behavior of some particular polynomials over $\mathbb{F}_{q^n}$. Next we try to connect the scattered behavior of $S(x) \in \mathbb{F}_{q^n}[x]$ of index $t(<r_1)$ with the permutation behavior of some particular polynomials over $\mathbb{F}_{q^n}$.
    \begin{theorem}
       Suppose $0<t<r_1<r_2< \cdots <r_k<n$, and $a_i \in \mathbb{F}_{q^n}^*$ such that order of $a_i$'s divide $(q^{t}-1), \forall i=1, \cdots, k$. Then $S(x)=\sum_{i=1}^{{k}} a_ix^{q^{r_i}} \in \mathbb{F}_{q^n}[x]$ is a SP of index t iff $S_{\rho}^t(x)= \sum_{i=1}^{{k}} a_i(\rho^{q^{r_i-t}}-\rho)x^{q^{r_i-t}} \in \mathbb{F}_{q^n}[x]$ is a permutation polynomial for any $ \rho \in \mathbb{F}_{q^n}^* - \mathbb{F}_q$.
   \end{theorem}
   \begin{proof}
       From Theorem(\ref{22}), we have $S(x)=\sum_{i=1}^{{k}} a_ix^{q^{r_i}} \in \mathbb{F}_{q^n}[x]$  is a SP of index $t(<r_1)$ \textit{iff} $S^t(x)=\sum_{i=1}^{{k}} a_ix^{q^{(r_i-t)}}$ is a SP of index 0. Also, using Proposition 2.6 of \cite{longobardi2021partially}, $S^t(x)=\sum_{i=1}^{{k}} a_ix^{q^{(r_i-t)}}$ is a SP of index 0 \textit{iff} $S_{\rho}^t(x)=S^t(\rho \cdot  x)-\rho \cdot S^t(x)  \in \mathbb{F}_{q^n}[x]$ is a permutation polynomial, for any $ \rho \in \mathbb{F}_{q^n}^*-\mathbb{F}_q$.\\
       So, $S_{\rho}^t(x)= \sum_{i=1}^{{k}} a_i(\rho^{q^{r_i-t}}-\rho)x^{q^{r_i-t}}$.\\
       Hence $S(x)=\sum_{i=1}^{{k}} a_ix^{q^{r_i}} \in \mathbb{F}_{q^n}[x]$ is a SP of index t \textit{iff} $S_{\rho}^t(x)= \sum_{i=1}^{{k}} a_i(\rho^{q^{r_i-t}}-\rho)x^{q^{r_i-t}} \in \mathbb{F}_{q^n}[x]$ is a permutation polynomial, for any $ \rho \in \mathbb{F}_{q^n}^*-\mathbb{F}_q$.
   \end{proof}
\section{Scattered Binomials}\label{sec4}
The binomials studied mostly in terms of their scattered behavior are the \textit{Lunardon-Polverino}  polynomials (known as LP polynomials) of the form $x^{q^r}+\delta x^{q^{n-r}} \in \mathbb{F}_{q^n}[x]$ such that $(n,r)=1, N_{q^n/q}(\delta) \neq 1$. Other binomials such as $x^{q^s}+\delta x^{q^{s(k-1)}}$ where $(s,k)=1$ and $N_{q^k/q}(\delta) \neq 1$ \cite{lunardon2001blocking}, $x^{q}+\delta x^{q^4} \in \mathbb{F}_{q^6}[x]$ where $q>4$ and for certain choice of $\delta$ \cite{BARTOLI2021111}, $x^{q}+\delta x^{q^5} \in \mathbb{F}_{q^8}[x]$ where $q$ is odd and $\delta^2=-1$ \cite{CSAJBOK2018203} are discussed in recent years. In this section, we focus on the linearized polynomials of the form $S(x)=a_1x^{q^{r_1}}+a_2x^{q^{r_2}} \in \mathbb{F}_{q^n}[x]$ such that $|a_2| \mid q^{r_1}-1$, and we relate the scattered behavior of such polynomial with the LP polynomial and the polynomials of the form $x^{q}+\delta x^{q^5} \in \mathbb{F}_{q^8}[x]$.
    \begin{proposition}\label{SB}
        Let $S(x)=a_1x^{q^{r_1}}+a_2x^{q^{r_2}} \in \mathbb{F}_{q^n}[x]$ such that $|a_2| \mid q^{r_1}-1$. Then the following are equivalent.
      \begin{enumerate}
          \item  $S(x)=a_1x^{q^{r_1}}+a_2x^{q^{r_2}} \in \mathbb{F}_{q^n}[x]$  is a SP of index $t \in \{r_1,r_2\}$,
          \item $(r_2-r_1,n)=1$.
      \end{enumerate}
    \end{proposition}
    \begin{proof}
         From Theorem (\ref{1.5}), $S(x)=a_1x^{q^{r_1}}+a_2x^{q^{r_2}}$ is a SP of index $r_1$ \textit{iff} $S_{r_1}^{r_1}(x)=a_2x^{q^{r_2-r_1}}$ is a SP of index 0 over $\mathbb{F}_{q^n}$. Further, $S_{r_1}^{r_1}(x)=a_2x^{q^{r_2-r_1}}$ is a SP of index 0 over $\mathbb{F}_{q^n}$ \textit{iff} $(n,r_2-r_1)=1$.\\
          Hence, $S(x)=a_1x^{q^{r_1}}+a_2x^{q^{r_2}}$ is a SP of index $r_1$ \textit{iff} $(n,r_2-r_1)=1$. \\
          Let $(y,z) \in C_i \times C_j \subset \mathbb{F}_{q^n}^* \times \mathbb{F}_{q^n}^*$. Then
        \begin{align*}
           \frac{S(y)}{y^{q^{r_2}}} = \frac{S(z)}{z^{q^{r_2}}}
           \Leftrightarrow \;  (\frac{y}{z})^{(q^{r_2}-q^{r_1})}=1
             \Leftrightarrow \; (\frac{y}{z})^{(q^{r_2-r_1}-1)}=1
           \Leftrightarrow \;  \frac{y^{q^{r_2-r_1}}}{y}=\frac{z^{q^{r_2-r_1}}}{z}\tag{9}\label{9}
        \end{align*}
        Suppose $S(x)$ is a SP of index $r_2$ with DP $(y_1,z_1)$. Then, from (\ref{1}), we have $$\frac{S(y_1)}{y_1^{q^{r_2}}} = \frac{S(z_1)}{z_1^{q^{r_2}}} \Rightarrow \frac{y_1}{z_1} \in \mathbb{F}_q.$$\\
        From (\ref{9}) we have $$\frac{y_1^{q^{r_2-r_1}}}{y_1}=\frac{z_1^{q^{r_2-r_1}}}{z_1} \Rightarrow \frac{y_1}{z_1} \in \mathbb{F}_q.$$\\
        So, $x^{q^{r_2-r_1}}$ is a SP of index 0 over $\mathbb{F}_{q^n}$, i.e., $(n,r_2-r_1)=1$.
        
        Conversely, let $(n,r_2-r_1)=1$. Then, $x^{q^{r_2-r_1}}$ is a SP of index 0 over $\mathbb{F}_{q^n}$.\\
        So, from (\ref{9}), $S(x)$ is a SP of index $r_2$ over $\mathbb{F}_{q^n}$.\\
        Hence, (1) and (2) are equivalent.
    \end{proof}
    \begin{example}
        Consider $S(x)=x^{25^{9}}+x^{25^{50}}\in \mathbb{F}_{25^{100}}[x]$. Then from Proposition (\ref{SB}), $S(x)$ is a SP of index $t \in \{9,50\}$.
    \end{example}
    \begin{example}
        Consider $S(x)=x^{101^{2}}+x^{101^{4}}\in \mathbb{F}_{101^{6}}[x]$. As $(4-2,6)=1$, then from Proposition (\ref{SB}), $S(x)$ is not a SP of index $t$ where $t\in \{2,4\}$.
    \end{example}
    Suppose $|a_1|,|a_2| \mid q^{r_1}-1$. In that case, we observe that for $t=r_2>r_1$, from Theorem (\ref{22})(2) and Proposition (\ref{SB}), $T(x)=a_1x+a_2x^{q^{r_2-r_1}}$ is a SP of index $\{r_2-r_1\}$ \textit{iff} $(n,r_2-r_1)=1$. Next, for $0<r<n$ and $a_1,a_2 \in \mathbb{F}_{q^n}^*$, we discuss the scattered behavior of $a_1x+a_2x^{q^r}$ explicitly over $\mathbb{F}_{q^n}$.
    \begin{theorem}\label{SB2}
        Let $r(<n) \in \mathbb{N}$ and  $a_1,a_2 \in \mathbb{F}_{q^n}^*$. Then the following are equivalent.
      \begin{enumerate}
          \item  $T(x)=a_1x+a_2x^{q^{r}} \in \mathbb{F}_{q^n}[x]$  is a SP of index $r$,
          \item $(r,n)=1$.          
      \end{enumerate}
    \end{theorem}
    \begin{proof}
        Let $(y,z) \in C_i \times C_j \subset \mathbb{F}_{q^n}^* \times \mathbb{F}_{q^n}^*$, then 
    \begin{align*}
          & \frac{T(y)}{y^{q^{r}}} = \frac{T(z)}{z^{q^{r}}}
           \Leftrightarrow \; \frac{a_1y+a_2y^{q^r}}{y^{q^r}}=\frac{a_1z+a_2z^{q^r}}{z^{q^r}}
           \Leftrightarrow \;  \frac{y^{q^r}}{y}=\frac{z^{q^r}}{z}\tag{10}\label{10}
        \end{align*}
Suppose $T(x)=a_1x+a_2x^{q^{r}} \in \mathbb{F}_{q^n}[x]$ is a SP of index $r$, and let $(y_1,z_1)$ be a DP in that case. Then from (\ref{1}), we have 
\begin{align*}
          & \frac{T(y_1)}{y_1^{q^{r}}} = \frac{T(z_1)}{z_1^{q^{r}}}
           \Rightarrow \frac{y_1}{z_1} \in \mathbb{F}_q.
        \end{align*}
        From (\ref{10}) we have
        \begin{align*}
          & \frac{y_1^{q^{r}}} {y_1}= \frac {z_1^{q^{r}}}{z_1}
           \Rightarrow \frac{y_1}{z_1} \in \mathbb{F}_q.
        \end{align*}
        So, $x^{q^r}$ is a SP of index 0, i.e., $(r,n)=1$.\\
        Conversely, let $(r,n)=1$. Then $x^{q^r}$ is a SP of index 0, i.e. there exists $(y_2,z_2) \in  \mathbb{F}_{q^n}^* \times \mathbb{F}_{q^n}^*$ such that  $\frac{y_2^{q^{r}}} {y_2}= \frac {z_2^{q^{r}}}{z_2}
           \Rightarrow \frac{y_2}{z_2} \in \mathbb{F}_q$. So, from (\ref{10}), we have $T(x)=a_1x+a_2x^{q^{r}} \in \mathbb{F}_{q^n}[x]$  is a SP of index $r$.\\
            Hence, (1) and (2) are equivalent.
    \end{proof}
    \begin{example}
        $x+x^{27^{81}}$ is a SP of index $81$ over $\mathbb{F}_{27^{110}}$, however, $x+x^{27^{80}}$ is not a SP of index $80$  over $\mathbb{F}_{27^{110}}$.
    \end{example}
    In this paper, we are considering a class of linearized polynomial of the form $S(x)=\sum_{i=1}^{{k}} a_ix^{q^{r_i}} \in \mathbb{F}_{q^n}[x]$ such that $1 
    \leq r_1<r_2< \cdots,<r_k<n$ and $a_i \in \mathbb{F}_{q^n}^*$, with $|a_i|\mid (q^{r_1}-1), \forall i=2, \cdots, k$. Next we work on the connection between the linearized polynomials of the form $S(x)=a_1x^{q^{r_1}}+a_2x^{q^{r_2}}$ and the \textit{Lunardon-Polverino}  polynomials (known as LP polynomials) of the form $x^{q^r}+\delta x^{q^{n-r}} \in \mathbb{F}_{q^n}[x]$ such that $(n,r)=1, N_{q^n/q}(\delta
    ) \neq 1$. \\
    Suppose $S_B(n,q)= \{ a_1x^{q^{r_1}}+a_2x^{q^{r}} \in \mathbb{F}_{q^n}[x]:a_i \in \mathbb{F}_{q^n}^*$, $\forall i=1,2$; with $|a_2|\mid (q^{r_1}-1), 1 \leq r_1<r_2<n\}$ and $S_{LP}(n,q)=\{x^{q^r}+\delta x^{q^{n-r}} \in \mathbb{F}_{q^n}[x]: \delta \in \mathbb{F}_{q^n}^*,(n,r)=1, N_{q^n/q}(\delta
    ) \neq 1 \}$. In general, we observe that $S_B(n,q),S_{LP}(n,q)$ are not directly related; however, with the help of few additional conditions, we get some relation of these two sets.
    \begin{lemma}\label{transf}
        Let $n(>1)$ be an odd integer such that $(q-1,n)=1$, and $\delta (\neq 1)\in \mathbb{F}_{q^n}^*$ with $|\delta|\mid (q-1)$. Then $S_B(n,q) \cap S_{LP}(n,q) \neq \phi$.
    \end{lemma}
    \begin{proof}
        For $a_1=1, a_2=\delta,r_1=r,r_2=n-r$ with $(n,r)=1$ we have $S(x)=x^{q^r}+\delta x^{q^{n-r}} \in S_B(n,q)$.\\
        If possible, let $N_{q^n/q}(\delta
    ) = 1$. Then from \cite{lidl1994introduction} we have $\delta=\phi^{q-1}$, for some $\phi \in \mathbb{F}_{q^n}^*$. \\
    So, $\delta^{(1+q+q^2+ \cdots+ q^{n-1})} =1$. As $|\delta|\mid (q-1)$, then we have $\delta^n=1$, i.e., $|\delta|\mid (q-1,n)$.\\
    As $(q-1,n)=1$, then we have $|\delta|\mid 1$, which is not possible. So, we must have $N_{q^n/q}(\delta
    ) \neq 1$. Now $S(x)=x^{q^r}+\delta x^{q^{n-r}} \in   \mathbb{F}_{q^n}[x]$ such that $(n,r)=1, N_{q^n/q}(\delta
    ) \neq 1$. So, $S(x)=x^{q^r}+\delta x^{q^{n-r}} \in S_{LP}(n,q)$.\\
    Hence $S_B(n,q) \cap S_{LP}(n,q) \neq \phi$.
    \end{proof}
    In general, if $S(x) \in S_B(n,q)$, then we do not know whether $S(x)$ is in $S_{LP}(n,q)$ or not. However, using the addition conditions in Lemma (\ref{transf}), we can transform $S(x)$ into an LP polynomial.
    From \cite{CSAJBOK2018203}, we know that for an odd $q$ and $\delta^2=-1$, $S(x)=x^q+\delta x^{q^5}$ is a scattered polynomial of index zero over $\mathbb{F}_{q^8}$. Further we study $S(x)$ and ask the follwing.\\
    (1) Can we discuss the scattered behavior of $S(x)$ for some positive index $t$ over $\mathbb{F}_{q^8}$?\\
    (2) Can we replace the condition $\delta^2=-1$ with any other condition (concerning $\delta$) to discuss the scattered behavior of $S(x)$ of index 0?\\
    In Theorem (\ref{csaj}), we observe that, with some additional condition on $q$, answer to the first question is yes. Later, we observe that, for some particular value of $q$, we can replace the condition $\delta^2=-1$.
    \begin{theorem}\label{csaj}
        Suppose $q$ is odd with $q \equiv1(\bmod \; {4})$ and $\delta (\neq -1,1) \in \mathbb{F}_{q^8}^*$ such that $\delta ^2=-1$. Then $S(x)=x^q+\delta x^{q^5}$ is not a scattered polynomial of index $t \in \{1,5\}$ over $\mathbb{F}_{q^8}$.
    \end{theorem}
 \begin{proof}
     Let $q-1=4k$, where $k \in \mathbb{N}$. As $\delta ^2=-1$, then $\delta^4=1$, i.e. $\delta^{4k}=\delta^{q-1}=1$.\\
     So, $|\delta| \mid (q-1)$. Now $S(x)=x^q+\delta x^{q^5} \in \mathbb{F}_{q^8}[x]$, where $|\delta| \mid q-1$. As $(5-1, 8)>1$, from Proposition (\ref{SB}), $S(x)=x^q+\delta x^{q^5} \in \mathbb{F}_{q^8}[x]$ is not a scattered polynomial of index $t \in \{1,5\}$.
 \end{proof}
 \begin{corollary}
     Let $\delta (\neq -1,1) \in \mathbb{F}_{5^8}^*$ such that $|\delta| \mid 4$, then $S(x)=x^5+\delta x^{5^5} \in \mathbb{F}_{5^8}[x]$ is a scattered polynomial of index 0.
 \end{corollary}
 \begin{proof}
     We know that $(\mathbb{F}_{5^8}^*, \cdot)$ is a cyclic group of even order. So, $(\mathbb{F}_{5^8}^*, \cdot)$ has $\phi(2)=1$ number of element of order 2, which is $-1$.\\
     As $|\delta| \mid 4$, then $|\delta^2| \mid 2$. If $|\delta^2|=1$, then $|\delta| \mid 2$ implies $\delta=-1$, which is not possible.\\
     So, $|\delta^2| = 2$. As $-1$ is the only element of order 2, we have $\delta^2=-1$. From \cite{CSAJBOK2018203}, for q=5, $S(x)=x^5+\delta x^{5^5} \in \mathbb{F}_{5^8}[x]$ is a scattered polynomial of index 0.     
 \end{proof}
   
   \section{A Family of Exceptional Scattered Binomial}
   Suppose $S(x)=\sum_{i=1}^{{k}} a_ix^{q^{r_i}} \in \mathbb{F}_{q^n}[x]$ is a SP of index $t \in \{0, \cdots, n-1\}$. Then $S(x)$ is said to be \textit{exceptional scattered} of index $t$ if there exist infinitely many $m \in \mathbb{N}$ such that, for any distinct $y,z \in  \mathbb{F}_{q^n}^*$, Condition (\ref{1}) holds. Recently in \cite{longobardi2021partially}, Longobardi and Zanella weakened the property of exceptional scatteredness of a polynomial by considering separately the
   exceptionality of two weaker properties, namely the L-$q^t$-partial scatteredness and the R-$q^t$-partial scatteredness. While several families of scattered polynomials have been constructed in recent years, below we present two well-known families of exceptional scattered polynomials.
   
   $\bullet$ $S(x)=x^{q^r}$ of index 0, with $(r,n)=1$ (polynomials of so-called pseudoregulus type);
   
   $\bullet$ $S(x)=x+\delta x^{q^{2r}}$ of index $r$, with $(r,n)=1$ and $N_{q^n/q}(\delta) \neq 1$ (so-called LP polynomials).\\
   In this section, using the exceptional scattered behavior of the LP polynomial of the form $x+\delta x^{q^{2r}}$, we present a family of exceptional scattered polynomial $x^q+\delta x^{q^{(2r+1)}}$ of index $\{r+1\}$. Next we show that, with some additional conditions, $x^q+\delta x^{q^{(2r+1)}}  \in \mathbb{F}_{q^n}[x]$ is a SP of index $\{r+1\}$.
\begin{theorem}
    Let $\delta (\neq 1) \in \mathbb{F}_{q^n}^*$ with $|\delta| \mid (q-1)$, and $n (>3)$ be odd such that $(n,q-1)=1$. If $(r,n)=1$ for some $r(<n) \in \mathbb{N}$, then $S(x)=x^q+\delta x^{q^{(2r+1)}}$ is a SP of index $\{r+1\}$ over $\mathbb{F}_{q^n}$. Further $S(x)=x^q+\delta x^{q^{(2r+1)}}$ is a SP of index $ t\in \{1,r+1,2r+1\}$ over $\mathbb{F}_{q^n}$.
\end{theorem}
\begin{proof}
    Let $T(x)=x+\delta x^{q^{2r}} \in \mathbb{F}_{q^n}[x]$. As $\delta (\neq 1) \in \mathbb{F}_{q^n}^*$, $n (>3)$ is odd, $(n,q-1)=1$, and $|\delta| \mid (q-1)$, from Lemma(\ref{transf}), we have  $N_{q^n/q}(\delta) \neq 1$. \\
    So, $N_{q^n/q}(\delta) \neq 1$ and $(r,n)=1$. From \cite{BARTOLI2022101956}, we have $T(x)=x+\delta x^{q^{2r}} \in \mathbb{F}_{q^n}[x]$ is a SP of index $\{r\}$. Let $r_1=1,r_2=2r+1,t=r+1$, then for $|\delta| \mid (q-1)$ we have $T(x)=x+\delta x^{q^{r_2-r_1}} \in \mathbb{F}_{q^n}[x]$ is a SP of index $\{t-r_1\}$. Then from Theorem(\ref{22}), we have $S(x)=x^{q^{r_1}}+\delta x^{q^{r_2}}$ is a SP of index $\{t\}$ over $\mathbb{F}_{q^n}$.\\
    Hence $S(x)=x^q+\delta x^{q^{(2r+1)}}$ is a SP of index $\{r+1\}$ over $\mathbb{F}_{q^n}$.

    For the second part, as $|\delta| \mid (q-1)$ and $(r,n)=1$, from Proposition(\ref{SB}), it follows that $S(x)=x^q+\delta x^{q^{(2r+1)}}$ is a SP of index $t\in \{1,r+1,2r+1\}$ over $\mathbb{F}_{q^n}$.
\end{proof}
\begin{theorem}\label{exp}
    Let $\delta (\neq 1) \in \mathbb{F}_{q^n}^*$ with $|\delta| \mid (q-1)$, and $n (>3)$ be odd such that $(n,q-1)=1=(r,n)$, for some $r(<n) \in \mathbb{N}$. Then $S(x)=x^q+\delta x^{q^{(2r+1)}}$ is an exceptional scattered polynomial of index $\{r+1\}$ over $\mathbb{F}_{q^n}$.
\end{theorem}
\begin{proof}
    Suppose $T(x)=x+\delta x^{q^{2r}} \in \mathbb{F}_{q^n}[x]$ and $m \in \mathbb{N}$. If $(y,z) \in  \mathbb{F}_{q^{mn}}^* \times \mathbb{F}_{q^{mn}}^*$ with $y \neq z$, then we have
        \begin{align*}
          & \frac{S(y)}{y^{q^{r+1}}} = \frac{S(z)}{z^{q^{r+1}}}\\
           \Leftrightarrow \; &  \frac{y^q+\delta^q y^{q^{(2r+1)}}}{y^{q^{r+1}}} = \frac{z^q+\delta^q z^{q^{(2r+1)}}}{z^{q^{r+1}}} [\text{as order of } \delta \text{ divides } (q-1) ]\\
           \Leftrightarrow \;  &  \frac{y+\delta y^{q^{(2r)}}}{y^{q^{r}}} = \frac{z+\delta z^{q^{(2r)}}}{z^{q^{r}}}\\
           \Leftrightarrow \; & \frac{T(y)}{y^{q^{r}}}=\frac{T(z)}{z^{q^{r}}}.\tag{11}\label{11}
        \end{align*} 
        As $\delta (\neq 1) \in \mathbb{F}_{q^n}^*$, $n (>3)$ is odd, $(n,q-1)=1$, and $|\delta| \mid (q-1)$, from Lemma(\ref{transf}), we have  $N_{q^n/q}(\delta) \neq 1$. Also, $(r,n)=1$.
    So, $T(x)=x+\delta x^{q^{2r}} \in \mathbb{F}_{q^n}[x]$ is an exceptional scattered polynomial of index $\{r\}$, i.e., there exists infinite number of $m \in \mathbb{N}$ with $(y_m,z_m) \in  \mathbb{F}_{q^{mn}}^* \times \mathbb{F}_{q^{mn}}^*$ and $y_m \neq z_m$, such that 
    \begin{align*}
        \frac{T(y_m)}{y_m^{q^{r}}} = \frac{T(z_m)}{z_m^{q^{r}}} \Rightarrow \frac{y_m}{z_m} \in \mathbb{F}_q.
    \end{align*}
     From (\ref{11}) we have 
    \begin{align*}
        \frac{S(y_m)}{y_m^{q^{r+1}}} = \frac{S(z_m)}{z_m^{q^{r+1}}} \Leftrightarrow \; \frac{T(y_m)}{y_m^{q^{r}}} = \frac{T(z_m)}{z_m^{q^{r}}} \Rightarrow \frac{y_m}{z_m} \in \mathbb{F}_q.
    \end{align*}
    So, there exists infinite number of $m \in \mathbb{N}$ with $(y_m,z_m) \in  \mathbb{F}_{q^{mn}}^* \times \mathbb{F}_{q^{mn}}^*$ and $y_m \neq z_m$, such that 
    \begin{align*}
        \frac{S(y_m)}{y_m^{q^{r+1}}} = \frac{S(z_m)}{z_m^{q^{r+1}}} \Rightarrow \frac{y_m}{z_m} \in \mathbb{F}_q.
    \end{align*}
    Hence $S(x)=x^q+\delta x^{q^{(2r+1)}}$ is an exceptional scattered polynomial of index $\{r+1\}$ over $\mathbb{F}_{q^n}$.
\end{proof}
    \bibliographystyle{plain}
    \bibliography{bibl.bib}
\end{document}